\documentclass[12pt,a4paper,english,intoc,bibliography=totoc,index=totoc,BCOR10mm,captions=tableheading,titlepage,fleqn]{scrbook}
\usepackage{lmodern}

\usepackage[T1]{fontenc}
\usepackage[latin9]{inputenc}
\usepackage{babel}
\usepackage{amsmath}
\usepackage{amssymb}
\usepackage[unicode=true,
 bookmarks=true,bookmarksnumbered=true,bookmarksopen=true,bookmarksopenlevel=2,
 breaklinks=false,pdfborder={0 0 0},backref=false,colorlinks=false]
 {hyperref}
\hypersetup{pdftitle={Abstract},
 pdfauthor={Uwe Stöhr},
 pdfpagelayout=OneColumn, pdfnewwindow=true, pdfstartview=XYZ, plainpages=false}
\usepackage{breakurl}

\makeatletter

\usepackage{calc}

\makeatother

\begin{document}
\noindent \begin{center}
\textbf{\LARGE A Borel maximal cofinitary group}
\par\end{center}{\LARGE \par}

\noindent \begin{center}
{\large Haim Horowitz and Saharon Shelah}
\par\end{center}{\large \par}

\noindent \begin{center}
\textbf{Abstract}
\par\end{center}

\noindent \begin{center}
We construct a Borel maximal cofinitary group.%
\footnote{Date: September 8, 2016.



Publication 10xx of the second author.

Partially supported by European Research Council grant 338821.%
}
\par\end{center}

\textbf{\large Introduction}{\large \par}

The study of mad families and their relatives occupies a central place
in modern set theory. As the straightforward way to construct such
families involves the axiom of choice, questions on the definability
of such families naturally arise. The following classical result is
due to Mathias:

\textbf{Theorem ({[}Ma{]}): }There are no analytic mad families.

In recent years, there has been considerable interest in the definability
of several relatives of mad families, such as maximal eventually different
families and maximal cofinitary groups. A family $\mathcal F \subseteq \omega^{\omega}$
is a maximal eventually different family if $f\neq g\in \mathcal F \rightarrow f(n) \neq g(n)$
for large enough $n$, and $\mathcal F$ is maximal with respect to
this property. The following result was recently discovered by the
authors:

\textbf{Theorem ({[}HwSh1089{]}): }Assuming $ZF$, there exists a
Borel maximal eventually different family.

As for maximal cofinitary groups (see definition 1 below), several
consistency results were established on the definability of such groups,
for example, the following results by Kastermans and by Fischer, Friedman
and Toernquist:

\textbf{Theorem ({[}Ka{]}): }There is a $\Pi^1_1$-maximal cofinitary
group in $L$.

\textbf{Theorem ({[}FFT{]}): }$\mathfrak b=\mathfrak c=\aleph_2$
is consistent with the existence of a maximal cofinitary group with
a $\Pi^1_2-$definable set of generators.

Our main goal in this paper is to establish the existence of a Borel
maximal cofinitary group in $ZF$. We intend to improve the current
results in a subsequent paper, and prove the existence of closed MED
families and MCGs.

$\\$

\textbf{\large The main theorem}{\large \par}

\textbf{Definition 1: }$G\subseteq S_{\infty}$ is a maximal cofinitary
group if $G$ is a subgroup of $S_{\infty}$, $|\{n : f(n)=n\}|<\aleph_0$
for every $Id \neq f \in G$, and $G$ is maximal with respect to
these properties.

\textbf{Theorem 2 $(ZF)$: }There exists a Borel maximal cofinitary
group.

The rest of the paper will be dedicated for the proof of the above
theorem. It will be enough to prove the existence of a Borel maximal
cofinitary group in $Sym(U)$ where $U$ is an arbitrary set of cardinality
$\aleph_0$.

Convention: Given two sequences $\eta$ and $\nu$, we write $\eta \leq \nu$
when $\eta$ is an initial segment of $\nu$.

\textbf{Definition 3: }The following objects will remain fixed throughout
the proof:

a. $T=2^{<\omega}$.

b. $\bar u=(u_{\rho} : \rho \in T)$ is a sequence of pairwise disjoint
sets such that $U=\cup \{u_{\rho} : \rho \in T\} \subseteq H(\aleph_0)$
(will be chosen in claim 4).

c. $<_*$ is a linear order of $H(\aleph_0)$ of order type $\omega$
such that given $\eta, \nu \in T$, $\eta<_* \nu$ iff $lg(\eta)<lg(\nu)$
or $lg(\eta)=lg(\nu) \wedge \eta <_{lex} \nu$.

d. For every $\eta \in T$, $\Sigma \{|u_{\nu}| : \nu <_* \eta\} \ll |u_{\eta}|$.

e. Borel functions $\bold B=\bold{B}_0$ and $\bold{B}_{-1}=\bold{B}_0^{-1}$
such that $\bold{B}: Sym(U) \rightarrow 2^{\omega}$ is injective
with a Borel image, and $\bold{B}_{-1}: 2^{\omega} \rightarrow Sym(U)$
satisfies $\bold{B}(f)=\eta \rightarrow \bold{B}_{-1}(\eta)=f$.

f. Let $\bold{A}_1=\{f\in Sym(U) : f$ has a finite number of fixed
points$\}$, $\bold{A}_1$ is obviously Borel.

g. $\{f_{\rho,\nu} : \nu \in 2^{lg(\rho)}\}$ generate the group $K_{\rho}$
(defined below) considered as a subgroup of $Sym(u_{\rho})$.

\textbf{Claim 4: }There exists a sequence $(u_{\rho},\bar{f_{\rho}}, \bar{A_{\rho}} : \rho \in T)$
such that:

a. $\bar{f_{\rho}}=(f_{\rho,\nu} : \nu \in T_{lg(\rho)})$.

b. $f_{\rho,\nu} \in Sym(u_{\rho})$ has no fixed points.

c. $\bar{A_{\rho}}=(A_{\rho,\nu} : \nu \in T_{lg(\rho)})$. We shall
denote $\underset{\nu \in T_{lg(\rho)}}{\cup}A_{\rho,\nu}$ by $A_{\rho}'$.

d. $A_{\rho,\nu} \subseteq u_{\rho} \subseteq H(\aleph_0)$ and $\Sigma \{|u_{\eta}| : \eta<_* \rho\} \ll |A_{\rho,\nu}|$.

e. $\nu_1 \neq \nu_2 \in T_{lg(\rho)} \rightarrow A_{\rho,\nu_1} \cap A_{\rho,\nu_2}=\emptyset$.

f. If $\rho \in 2^n$ and $w=w(...,x_{\nu},...)_{\nu \in 2^n}$ is
a non-trivial group term of length $\leq n$ then:

1. $w(...,f_{\rho,\nu},...)_{\nu \in 2^n} \in Sym(u_{\rho})$ has
no fixed points.

2. $(w(...,f_{\rho,\nu},...)_{\nu \in 2^n}(\underset{\nu \in 2^n}{\cup}A_{\rho,\nu})) \cap \underset{\nu \in 2^n}{\cup}A_{\rho,\nu}=\emptyset$.

g. $\{f_{\rho,\nu} : \nu \in T_{lg(\rho)}\}$ generate the group $K_{\rho}$
(whose set of elements is $u_{\rho}$) which is considered as a group
of permutations of $u_{\rho}$.

\textbf{Proof: }We choose $(u_{\rho}, \bar{f_{\rho}},\bar{A_{\rho}})$
by $<_*$-induction on $\rho$ as follows: Arriving at $\rho$, we
choose the following objects:

a. $n_{\rho}^1$ such that $\Sigma \{|u_{\nu}| : \nu <_* \rho\}2^{lg(\rho)+7} \ll n_{\rho}^1$
and let $n_{\rho}^0=\frac{n_{\rho}^1}{2^{lg(\rho)}}$.

b. Let $H_{\rho}$ be the group generated freely by $\{x_{\rho,\nu} : \nu \in T_{lg(\rho)}\}$.

c. In $H_{\rho}$ we can find $(y_{\rho,n} : n<\omega)$ which freely
generate a subgroup (we can do it explicitly, for example, if $a$
and $b$ freely generate a group, then $(a^nb^n : n<\omega)$ are
as required), wlog for $w_1$ and $w_2$ as in 4(f) and $n_1<n_2$
we have $w_1y_{\rho,n_1} \neq w_2y_{\rho,n_2}$. 

Now choose $A_{\rho,\nu}^1 \subseteq \{y_{\rho,n} : n<\omega\}$ for
$\nu \in 2^{lg(\rho)}$ such that $\nu_1 \neq \nu_2 \rightarrow A_{\rho,\nu_1}^1 \cap A_{\rho,\nu_2}^1=\emptyset$
and $n_{\rho}^0 \leq |A_{\rho,\nu}^1|$.

d. Let $\Lambda_{\rho}=\{w : w=w(...,x_{\rho,\nu},...)_{\nu \in T_{lg(\rho)}}$
is a group word of length $\leq lg(\rho) \}$.

As $H_{\rho}$ is free, it's residually finite, hence there is a finite
group $K_{\rho}$ and an epimorphism $\phi_{\rho}: H_{\rho} \rightarrow K_{\rho}$
such that $\phi_{\rho} \restriction ((\underset{\nu \in 2^{lg(\rho)}}{\cup}A_{\rho,\nu}^1) \cup \Lambda_{\rho} \cup \{wa : w\in \Lambda_{\rho} \wedge a\in \underset{\nu \in 2^{lg(\rho)}}{\cup}A_{\rho,\nu}^1\})$
is injective (note that there is no use of the axiom of choice as
we can argue in a model of the form $L[A]$). WLOG $K_{\rho} \subseteq H(\aleph_0)$
and $K_{\rho}$ is disjoint to $\cup \{u_{\nu} : \nu <_* \rho\}$. 

We now define the following objects:

a. $u_{\rho}=K_{\rho}$.

b. $A_{\rho,\nu}=\{\phi_{\rho}(a) : a \in A_{\rho,\nu}^1\}$.

c. For $\nu \in 2^{lg(\rho)}$, let $f_{\rho,\nu} : u_{\rho} \rightarrow u_{\rho}$
be multiplication by $\phi_{\rho}(x_{\rho,\nu})$ from the left.

It's now easy to verify that $(u_{\rho}, \bar{A_{\rho}}, \bar{f_{\rho}})$
are as required, so $U=\cup \{u_{\rho} : \rho \in T\}$. $\square$

\textbf{Definition and claim 5: }A. a. Given $f\in Sym(U)$, let $g=F_1(f)$
be $g_{\bold{B}(f)}^*$, where for $\nu \in 2^{\omega}$, $g_{\nu}^*$
is the permutation of $U$ defined by: $g^*_{\nu}\restriction u_{\rho}=f_{\rho,\nu \restriction lg(\rho)}$
(recall that $\bar u$ is a partition of $U$ and each $f_{\rho,\nu}$
belongs to $Sym(u_{\rho})$, therefore $g$ is well-defined and belongs
to $Sym(U)$).

b. Let $G_1$ be the subgroup of $Sym(U)$ generated by $\{g_{\nu}^* : \nu \in 2^{\omega}\}$
(which includes $\{F_1(f) : f\in Sym(U)\}$).

c. Let $I_1$ be the ideal on $U$ generated by the sets $v\subseteq U$
satisfying the following property:

$(*)_v$ For some $\rho=\rho_v \in 2^{\omega}$, for every $n$, there
is at most one pair $(a,\nu)$ such that $\nu \in T$, $a\in v\cap u_{\nu}$
and $\rho \cap \nu=\rho \restriction n$. 

c(1). Note that $I_1$ is indeed a proper ideal: Suppose that $v_0,...,v_n$
are as above and let $\rho_0,...,\rho_n$ witness $(*)_{v_i}$ $(i=0,...,n)$.
Choose $k$ such that $2^k>n+1$ and choose $\eta \in 2^k \setminus \{\rho _i \restriction k : i\leq n\}$.
For each $i\leq n$, there is $k(i) \leq k$ such that $\eta \cap \rho_i= \rho_i \restriction k(i)$.
For each $i\leq n$, let $n(i)$ be the length $\nu$ such that $(a,\nu)$
witness $(*)_{v_i}$ for $k(i)$. Choose $\eta'$ above $\eta$ such
that $lg(\eta')>n(i)$ for every $i$, then $u_{\eta'} \cap (\underset{i\leq n}{\cup}v_i)=\emptyset$.

d. Let $K_1=\{f \in Sym(U) : fix(f) \in I_1 \}$ where $fix(f)=\{x : f(x)=x\}$.

e. For $\eta \in T$, $a,b\in u_{\eta}$, $n=lg(\eta)<\omega$ and
let $y_{a,b}=((f_{\eta,\rho_{a,b,l}},i_{a,b,l}) : l<l_{a,b}=l(*))$
such that:

1. $\rho_{a,b,l} \in 2^n$.

2. $i_{a,b,l} \in \{1,-1\}$.

3. $b=(f_{\eta,\rho_{a,b,0}})^{i_{a,b,0}} \cdot \cdot \cdot (f_{\eta,\rho_{a,b,l(*)-1}})^{i_{a,b,l(*)-1}}(a)$.

4. $l_{a,b}=l(*)$ is minimal under 1-3, $y_{a,b}$ is $<_*$-minimal
under this requirement.

5. $i_l \neq i_{l+1} \rightarrow \rho_{a,b,l} \neq 
\rho_{a,b,l+1}$.

By claim 4(g) and definition 3(c), $y_{a,b}$ is always well-defined.

B. There are Borel functions $\bold{B}_{1,1}, \bold{B}_{1,2},$ etc
with domain $Sym(U)$ such that:

a. $\bold{B}_{1,1}(f) \in \{0,1\}$ and $\bold{B}_{1,1}(f)=0$ iff
$|fix(f)|<\aleph_0$.

b. Letting $\eta_1=\bold{B}(f)$, $\bold{B}_{1,2}(f) \in \{0,1\}$
and $\bold{B}_{1,2}(f)=1$ iff $\bold{B}_{1,1}(f)=0$ and for infinitely
many $n$, $f''(A_{\eta_1 \restriction n}') \nsubseteq \cup \{u_{\rho} : \rho \leq_* \eta_1 \restriction n\}$
(where $A_{\eta_1 \restriction n}'$ is defined in 4(c).

c. $\bold{B}_{1,3}(f) \in \omega$ such that if $\bold{B}_{1,1}(f)=\bold{B}_{1,2}(f)=0$
then for every $\bold{B}_{1,3}(f) \leq n$, $f''(A_{\eta_1 \restriction n}') \subseteq \cup \{u_{\rho} : \rho \leq_* \eta_1 \restriction n\}$.

d. $\bold{B}_{1,4}(f) \in \{0,1\}$ and $\bold{B}_{1,4}(f)=1$ iff
$\bold{B}_{1,1}(f)=\bold{B}_{1,2}(f)=0$ and $\{l_{a,f(a)} : a\in v_n$
and $\bold{B}_{1,3}(f) \leq n\}$ is unbounded, where $v_n:=\{a\in A_{\eta_1 \restriction n}' \subseteq u_{\eta_1 \restriction n} : f(a) \in u_{\eta_1 \restriction n}\}$.

e. $\bold{B}_{1,5}(f) \in \omega$ such that: If $\bold{B}_{1,4}(f)=\bold{B}_{1,2}(f)=\bold{B}_{1,1}(f)=0$
then $\bold{B}_{1,5}(f)$ is a bound of $\{l_{a,f(a)} : a\in v_n$
and $\bold{B}_{1,3}(f) \leq n\}$.

f. $\bold{B}_{1,6}(f) \in \{0,1\}$ such that: $\bold{B}_{1,6}(f)=1$
iff $\bold{B}_{1,1}(f)=\bold{B}_{1,2}(f)=\bold{B}_{1,4}(f)=0$ and
for every $m$ there exists $n>m$ such that: There are $a_1 \neq a_2 \in v_n$
such that for some $l$, $l<min\{l_{a_1,f(a_1)}, l_{a_2,f(a_2)}\}$,
$\rho_{a_1,f(a_1),l} \neq \rho_{a_2,f(a_2),l}$ and $\rho_{a_1,f(a_1),l} \restriction m=\rho_{a_2,f(a_2),l} \restriction m$.

g. $\bold{B}_{1,7}(f)$ is a sequence $(a_n=a_n(f) : n\in \bold{B}_{1,8}(f))$
such that if $\bold{B}_{1,6}(f)=1$ then:

1. $a_n \in v_n$

2. $\bold{B}_{1,8}(f) \in [\omega]^{\omega}$

3. $l_{a_n,f(a_n)}=l(*)=\bold{B}_{1,9}(f)$

4. $l_{**}=\bold{B}_{2,0}(f)<l_*$

5. $(\rho_{a_n,f(a_n),l_{**}} : n\in \bold{B}_{1,8}(f))$ are pairwise
incomparable.

6. For every $l<l_*$, the following sequence is constant: $(TV(\rho_{a_n,f(a_n),l} \leq \rho_{a_k,f(a_k),l}) : n<k \in \bold{B}_{1,8}(f))$.

h. $\bold{B}_{2,1}(f)$ is a sequence $(A_n=A_n(f) : n\in \bold{B}_{2,2}(f))$
such that if $\bold{B}_{1,1}(f)=\bold{B}_{1,2}(f)=\bold{B}_{1,4}(f)=\bold{B}_{1,6}(f)=0$
then:

1. $\bold{B}_{2,2}(f) \in [\omega]^{\omega}$

2. $A_n \subseteq A_{\eta_1 \restriction n}'$ (recalling that $\eta_1=\bold{B}(f)$)

3. $l_{a,f(a)}=l_*=\bold{B}_{2,3}(f)$ for $n\in \bold{B}_{2,2}(f)$
and $a\in A_n$

4. $(i_l : l<l_*)=(i_{a,f(a),l} : l<l_*)$ (recalling definition 5(e))
for every $n\in \bold{B}_{2,2}(f)$ and $a\in A_n$.

5. $\frac{1}{\bold{B}_{2,3}'(f)(n)} \leq \frac{|A_n|}{|v_n|}$ where
$\bold{B}_{2,3}'(f)(n) \in \omega \setminus \{0\}$, $\bold{B}_{2,3}(f)(n) \ll |v_n|$
and $v_n$ is defined in 5(B)(d).

6. $\bold{B}_{2,4,n}(f)=\bar{\rho_n^*}=(\rho_l^n : l<l_*)=(\rho_{a,f(a),l} : l<l_*)$
for every $n\in \bold{B}_{2,2}(f)$ and $a\in A_n$.

7. $(TV(\rho_l^n \leq \rho_l^m) : n<m \in \bold{B}_{2,2}(f))$ is
constantly $\bold{B}_{2,5,l}(f)$

i. $\bold{B}_{2,6}(f) \in \{0,1\}$ is $1$ iff $\bold{B}_{1,1}(f)=\bold{B}_{1,2}(f)=\bold{B}_{1,4}(f)=\bold{B}_{1,6}(f)=0$
and in $(h)(7)$, $\bold{B}_{2,5,l}=false$ for some $l<l_*$.

j. $\bold{B}_{2,6}'(f) \in \{0,1\}$ is $0$ iff $\bold{B}_{1,1}(f)=\bold{B}_{1,2}(f)=\bold{B}_{1,4}(f)=\bold{B}_{1,6}(f)=0$
and $\bold{B}_{2,6}(f)=0$

\textbf{Proof: }By the proof of Ramsey's theorem and the arguments
which are implicit in the proof of claim 7 below. Note that while
the statement {}``there exists an infinite homogeneous set'' is
analytic, we can Borel-compute that homogeneous set. See the proof
of claim 6 in {[}HwSh:1089{]} for more details. $\square$

\textbf{Definition and claim 6: }a. 1.\textbf{ }Let $H_3$ be the
set of $f\in Sym(U)$ such that:

$\alpha$. $\bold{B}_{1,1}(f)=0$

$\beta$. If $\bold{B}_{1,2}(f)=\bold{B}_{1,4}(f)=\bold{B}_{1,6}(f)=0$
then $\bold{B}_{2,6}'(f)=1$

a. 2. $H_3$ is Borel.

b. For $f\in Sym(U)$ let $G_f$ be the set of $g\in Sym(U)$ such
that:

1. If $f\notin H_3$ then $G_f=\{F_1(f)\}$.

2. If $f\in H_3$ then $G_f$ be the set of $g$ such that for some
$(B,\eta_1,\eta_2,\bar a, \bar b, \bar c, \bar d, \bar e,\bar{\nu})$
we have:

A. $B\subseteq \omega$ is infinite.

B. $\eta_1=\bold{B}(f) \in 2^{\omega}$ and $\eta_2 \in 2^{\omega}$.

C. $\bar a=(a_n : n\in B)$.

D. If $n\in B$ then $a_n \in A_{\eta_1 \restriction n}'=\underset{\rho \in 2^n}{\cup}A_{\eta_1 \restriction n,\rho} \subseteq u_{\eta_1 \restriction n}$
(recall that we denote $\underset{\rho \in 2^n}{\cup}A_{\eta_1 \restriction n,\rho}$
by $A_{\eta_1 \restriction n}'$).

E. $\bar b=(b_n : n\in B)$ and $\bar{\nu}=(\nu_n : n\in B)$, $\nu_n \in T$,
such that for each $n\in B$, $b_n=f(a_n)$ and $b_n \in u_{\nu_n}$.
$\bar c=(c_n : n \in B)$, $\bar d=(d_n : n \in B)$ and $\bar e=(e_n : n\in B)$
are such that $b_n,c_n \in u_{\nu_n}$ and $e_n \in u_{\eta_1 \restriction n}$.

F. For every $n\in B$, $g(a_n)=f(a_n)=b_n$.

G. For every $n\in B$, $g(b_n)=F_1(f)(a_n)=e_n$.

H. For every $n\in B$ we have $c_n=F_1(f)^{-1}(f(a_n))$ and $g(c_n)=F_1(f)(f(a_n))=d_n$.

I. If $b \in U$ is not covered by clauses F-H, then $g(b)=F_1(f)(b)$.

J. $g$ has no fixed points.

K. One of the following holds:

a. For every $n\in B$, $\eta_1 \restriction n <_* \nu_n$, $lg(\eta_2 \cap \nu_n)>max\{lg(\nu_m) : m\in B\cap n\}$
hence $(\nu_n \cap \eta_2 : n\in B)$ is $\leq-$increasing.

b. For every $n\in B$, $\nu_n=\eta_1 \restriction n$ and $l(a_n,f(a_n),n)$
is increasing (see definition 5(e)).

c. For every $n\in B$, $\nu_n=\eta_1 \restriction n$ and in addition,
$l(a_n,f(a_n),n)=l_*$ for every $n$, $i_{a_n,f(a_n),l}=i_l$ for
$l<l_*$ and for some $l_{**}<l_*$, the elemnts of $(\rho_{a_n,f(a_n),l_{**}} : l_{**}<l_*)$
are pairwise incomparable. $\square$

\textbf{Claim 7: }If $f\in H_3$ then there exists $g\in Sym(U)$
such that for some $(B,\eta_1,\eta_2,\bar a,\bar b,\bar{\nu})$, $g$
and $(B,\eta_1,\eta_2,\bar a,\bar b, \bar{\nu})$ are as required
in claim 6(c)(2) (and therefore, there are also $(\bar c,\bar d,\bar e)$
as required there). Moreover, $g$ is unique once $(B,\eta_1,\eta_2,\bar a,\bar b,\bar{\nu})$
is fixed.

\textbf{Remark 7A: }In claim 9 we need $g$ to be Borel-computable
from $f$, which is indeed the case by the discussion in the proof
of claim 5 and by the proof of claim 6 in {[}HwSh:1089{]}.

\textbf{Proof: }$f\in H_3$, so $\bold{B}_{1,1}(f)=0$.

We shall first observe that if $g$ is defined as above, then $g$
is a permutation of $U$ with no fixed points. It's also easy to see
that $g$ is unique once $(B,\eta_1,\eta_2,\bar a,\bar b,\bar{\nu})$
has been chosen. Therefore, it's enough to find $(B,\eta_1,\eta_2,\bar a,\bar b,\bar{\nu})$
as required.

\textbf{Case I ($\bold{B}_{1,2}(f)=1$): For infinitely many $n$,
$f''(A_{\eta_1 \restriction n}') \nsubseteq \cup \{u_{\rho} : \rho \leq_* \eta_1 \restriction n\}$.
}In this case, let $B_0=\{n : $ there is $a\in A_{\eta_1 \restriction n}'$
such that $f(a) \notin \cup \{u_{\rho} : \rho \leq_* \eta_1 \restriction n\}$,
and for every $n\in B_0$, let $a_n$ be the $<_*$-first element
in $A_{\eta_1 \restriction n}'$ witnessing that $n\in B_0$. Let
$b_n=f(a_n)$ and let $\nu_n \in T$ be the sequence for which $b_n \in u_{\nu_n}$.
Apply Ramsey's theorem (we don't need the axiom of choice, as we can
argue in some $L[A]$) to get an infinite set $B \subseteq B_0$ such
that $c_{k,l} \restriction [B]^k$ is constant for every $(k,l) \in \{(2,1),(2,2),(2,4),(3,1),(3,3)\}$,
where for $n_1<n_2<n_3$:

a) $c_{2,1}(n_1,n_2)=TV(lg(\nu_1)<lg(\nu_2))$.

b) $c_{2,2}(n_1,n_2)=TV(\nu_{n_2} \in \{\nu_n : n\leq n_1\})$.

c) $c_{3,1}(n_1,n_2,n_3)=TV(lg(\nu_{n_2} \cap \nu_{n_3})>\nu_{n_1})$.

d) $c_{3,3}(n_1,n_2)=\nu_{n_2}(lg(\nu_{n_1} \cap \nu_{n_2})) \in \{0,1,$undefined$\}$.

We shall prove now that $(lg(\nu_n) : n\in B)$ has an infinite increasing
subsequence: Choose an increasing sequence $n(i) \in B$ by induction
on $i$ such that $j<i \rightarrow lg(\nu_{n(j)})<lg(\nu_{n(i)})$.
Arriving at stage $i=j+1$, suppose that there is no such $n(i)$,
then $\{f(a_n) : n\in B \setminus n(j)\} \subseteq \cup \{u_{\rho} : lg(\rho) \leq lg(\nu_{n(j)})\}$,
hence $\{f(a_n) : n\in B\setminus n(j)\}$ is finite. Similarly, $\{\nu_n : n\in B\setminus n(j)\}$
is finite, and therefore, there are $n_1<n_2 \in B\setminus n(j)$
such that $\nu_{n_1}=\nu_{n_2}$ and $f(a_{n_1})=f(a_{n_2})$. As
$f$ is injective, $a_{n_1}=a_{n_2}$, and by the choice of the $a_n$,
$a_{n_1} \in u_{\eta_1 \restriction n_1}$ and $a_{n_2} \in u_{\eta_1 \restriction n_2}$.

This is a contradiction, as $u_{\eta_1 \restriction n_1} \cap u_{\eta_1 \restriction n_2}=\emptyset$.

Therefore, there is an infinite $B' \subseteq B$ such that $(lg(\nu_n) : n\in B')$
is increasing, and wlog $B'=B$.

Now we shall note that if $n_1<n_2<n_3$ are from $B$, then $lg(\nu_{n_2} \cap \nu_{n_3})>lg(\nu_{n_1})$: 

By the choice of $B$, $c_{3,1}(n_1,n_2,n_3)$ is constant for $n_1<n_2<n_3$,
so it suffices to show that $c_{3,1} \restriction [B]^3=true$. Let
$n_1=min(B)$ and $k=lg(\nu_{n_1})+1$. The sequence $(\nu_n \restriction k : n\in B\setminus \{n_1\})$
is infinite, hence there are $n_2<n_3 \in B\setminus \{n_1\}$ such
that $\nu_{n_2} \restriction k=\nu_{n_3} \restriction k$. Therefore,
$lg(\nu_{n_1})<k \leq lg(\nu_{n_2} \cap \nu_{n_3})$, and as $c_{3,1}$
is constant on $[B]^3$, we're done.

For $n<k \in B$ such that $k$ is the successor of $n$ in $B$,
let $\eta_n=\nu_n \cap \nu_k$. Suppose now that $n<k<l$ are successor
elements in $B$, then $lg(\eta_k)=lg(\nu_k \cap \nu_l)>lg(\nu_n) \geq lg(\nu_n \cap \nu_k)=lg(\eta_n)$,
and $\eta_n,\eta_k \leq \nu_k$, therefore, $\eta_n$ is a proper
initial segment of $\eta_k$ and $\eta_2:=\underset{n<\omega}{\cup}\eta_n \in 2^{\omega}$.
If $n<k \in B$ are successor elements, then $lg(\eta_k)>lg(\nu_n)$
(by a previous claim), therefore, $\nu_n \cap (T\setminus \nu_k)$
is disjoint to $\eta_2$, hence $\nu_n \cap \eta_2=\nu_n \cap \nu_k=\eta_n$.
Therefore, if $n<k \in B$ then $\nu_n \cap \eta_2$ is a proper initial
segment of $\nu_k \cap \eta_2$ and $\eta_2=\underset{n<\omega}{\cup}(\nu_n \cap \eta_2)$.

It's now easy to verify that $(B,\eta_1,\eta_2,\bar a,\bar b,\bar{\nu})$
and $g$ are as required.

\textbf{Case II ($\bold{B}_{1,2}(f)=0$ and $n_1$ stands for $\bold{B}_{1,3}(f)$):
There is $n_1$ such that for every $n_1 \leq n$, $f''(A_{\eta_1 \restriction n}') \subseteq \cup \{u_{\rho} : \rho \leq_* \eta_1 \restriction n\}$. }

For each $n$, recall that $v_n=\{a\in A_{\eta_1 \restriction n}' \subseteq u_{\eta_1 \restriction n} : f(a)\in u_{\eta_1 \restriction n}\}$.
$v_n$ satisfies $|A_{\eta_1 \restriction n}' \setminus v_n| \leq \Sigma \{|u_{\nu}| : \nu <_*\ \eta_1 \restriction n\}$,
and as $\Sigma\{|u_{\nu}| : \nu<_* \eta_1 \restriction n\} \ll |A_{\eta_1 \restriction n}'|$,
it follows that $\Sigma\{|u_{\nu}| : \nu<_* \eta_1 \restriction n\} \ll |v_n|$.
Recall also that for $a\in v_n$, as $f(a) \in u_{\eta_1 \restriction n}$,
by definition 5(e), $y_{a,f(a)}$ is well-defined.

We now consider three subcases:

\textbf{Case IIA ($\bold{B}_{1,4}(f)=1$): The set of $l_{a,f(a),n}$
for $a\in v_n$ and $n_1 \leq n$ is unbounded. }In this case, we
find an infinite $B\subseteq [n_1,\omega)$ and $a_n \in v_n$ for
each $n\in B$ such that $(l_{a_n,f(a_n),n} : n\in B)$ is increasing.
Now let $\eta_2:=\eta_1$ and define $\bar b$, $\bar{\nu}$ and $g$
as described in Definition 6. It's easy to verify that $(B,\eta_1,\eta_2,\bar a,\bar b, \bar{\nu})$
are as required.

\textbf{Case IIB ($\bold{B}_{1,4}(f)=0$ and $\bold{B}_{1,6}(f)=1$):
Case IIA doesn't hold, but $\bold{B}_{1,6}(f)=1$ and there is an
infinite $B\subseteq [n_1,\omega)$, $l_{**}<l_*$ (see below) and
$(a_n \in v_n : n\in B)$ (given by $\bold{B}_{1,8}(f)$, $\bold{B}_{2,0}(f)$
and $\bold{B}_{1,7}(f)$, respectively) such that:}

\textbf{a. $l_{a_n,f(a_n)}=l_*$ and \textbf{$l_{**}=\bold{B}_{2,0}(f)<l_*$}.}

\textbf{b. $i_{a_n,f(a_n),l}=i_l^*$ for $l<l_*$.}

\textbf{c. If $n\in B$ and $m\in B\cap n$, then $\rho_{a_m,f(a_m),l_{**}} \nsubseteq \rho_{a_n,f(a_n),l_{**}}$.}

In this case we define $\bar b$, $\bar{\nu}$ and $g$ as in Definition
6 and we let $\eta_2:=\eta_1$. It's easy to see that $(B,\eta_1,\eta_2,\bar a,\bar b \bar{\nu})$
are as required.

Remark: By a routine Ramsey-type argument, it's easy to prove that
if $\bold{B}_{1,6}(f)=1$ then the values of $\bold{B}_{1,7}(f),\bold{B}_{2,0}(f)$
are well-defined and Borel-computable so the above conitions hold.

\textbf{Case IIC ($\bold{B}_{1,4}(f)=\bold{B}_{1,6}(f)=0$): $\neg IIA \wedge \neg IIB$.
}We shall first prove that $\bold{B}_{2,1}(f)$, $\bold{B}_{2,2}(f)$,
$\bold{B}_{2,3}(f)$, $(\bold{B}_{2,4,n}(f) : n\in \bold{B}_{2,2}(f))$
and $(\bold{B}_{2,5,l}(f) : l<\bold{B}_{2,3}(f))$ are well-define
and Borel computable.

Let $l(*)$ be the supremum of the $l(a,f(a))$ where $n_1 \leq n$
and $a\in v_n$ ($l(*)<\omega$ by $\neg2A$). We can find $l(**) \leq l(*)$
such that $B_1:=\{ n\in B : v_{n,1}=\{a\in v_n : l(a,f(a))=l(**)\}$
has at least $\frac{v_n}{l(*)}$ elements$\}$ is infinite. Next,
we can find $i_*(l) \in \{1,-1\}$ for $l<l(**)$ such that $B_2:\{n \in B_1 : v_{n,2}=\{a\in v_{n,1} : \underset{l<l(**)}{\wedge}i_{a,f(a),l}=i_*(l)\}$
has at least $\frac{|v_n|}{2^{l(**)}l(*)}$ elements$\}$ is infinite.
For each $n\in B_2$, there are $\rho_{n,0},...,\rho_{n,l(**)-1} \in T_n$
such that $v_{n,3}=\{a\in v_{n,2} : \underset{l<l(**)}{\wedge}\rho_{a,f(a),l}=\rho_{n,l} \}$
has at least $\frac{|v_n|}{l(**)2^{l(**)}2^{nl(**)}}$ elements. By
Ramsey's theorem, there is an infinite subset $B_3 \subseteq B_2$
such that for each $l<l(**)$, the sequence $(TV(\rho_{m,l} \leq \rho_{n,l}) : m<n \in B_3)$
is constant. Therefore, we're done showing that the above Borel functions
are well-defined.

Now if $\bold{B}_{2,6}'(f)=1$ then we finish as in the previous case
(this time we're in the situation of 6(b)(2)(K)(c)). If $\bold{B}_{2,6}'(f)=0$,
then we get a contradiction to the assumption that $f\in H_3$, therefore
we're done. $\square$

\textbf{Claim 8: }If $w(x_0,...,x_{k_*-1})$ is a reduced non-trivial
group word, $f_0,...,f_{k_{**}-1} \in H_3$ are pairwise distinct,
$g_l \in G_{f_l}$ $(l\in \{0,...,k_{**}-1\})$, $g_l=g_{\nu_l}^*$
where $\nu_l \in 2^{\omega} \setminus \{\bold{B}(f) : f\in H_3\}$,
$l=k_{**},...,k_*-1$ and $(\nu_l : k_{**} \leq l<k_* )$ is without
repetition, then $w(g_0,...g_{k_*-1}) \in Sym(U)$ has a finite number
of fixed points.

Notation: For $l<k_{**}$, let $\nu_l:=\bold{B}(f_l)$.

\textbf{Proof: }Assume towards contradiction that $w(x_0,...,x_{k_*-1})=x_{k(m-1)}^{i(m-1)} \cdot...\cdot x_{k(0)}^{i(0)}$,
$\{f_0,...,f_{k_*-1}\}$ and $\{g_0,...,g_{k_*-1}\}$ form a counterexample,
where $i(l) \in \{-1,1\}$, $k(l)<k_*$ and $k(l)=k(l+1) \rightarrow \neg (i(l)=-i(l+1))$
for every $l<m$. WLOG $m=lg(w)$ is minimal among the various countrexamples.
Let $C=\{a\in U : w(g_0,...,g_{k_*-1})(a)=a\}$, this set is infinite
by our present assumption. For $c\in C$, define $b_{c,l}$ by induction
on $l<m$ as follows:

1. $b_{c,0}=c$.

2. $b_{c,l+1}=g_{k(l)}^{i(l)}(b_{c,l})$.

Notational warning: The letter $c$ with additional indices will be
used to denote the elements of sequences of the form $\bar c$ from
Definition 6(b)(2).

For all but finitely many $c\in C$, $(b_{c,l} : l<m)$ is without
repetition by the minimality of $m$, so wlog this is true for every
$c\in C$.

For every $c\in C$, let $\rho_{c,l} \in T$ be such that $b_{c,l} \in u_{\rho_{c,l}}$,
and let $l_1[c]$ be such that $\rho_{c,l_1[c]} \leq_* \rho_{c,l}$
for every $l<m$. We can choose $l_1[c]$ such that one of the following
holds:

1. $l_1[c]>0$ and $\rho_{c,l_1[c]-1} \neq \rho_{c,l_1[c]}$

2. $l_1[c]=0$ and $\rho_{c,m-1} \neq \rho_{c,0}$

3. $\rho_{c,0}=...=\rho_{c,m-1}$

We may assume wlog that $(l_1[c] : c\in C)$ is constant and that
actually $l_1[c]=0$ for every $c\in C$. In order to see that we
can assume the second part, for $j<m$ let $w_j(x_0,...,x_{k_*-1})=x_{k(j-1)}^{i(j-1)} \cdot \cdot \cdot x_{k(0)}^{i(0)}x_{k(m-1)}^{i(m-1)} \cdot \cdot \cdot x_{k(j)}^{i(j)}$,
then $w_j(g_0,...,g_{k_*-1}) \in Sym(U)$ is a conjugate of $w(g_0,...,g_{k_*-1})$.
The set of fixed points of $w_j(g_0,...,g_{k_*-1})$ includes $\{b_{c,j} : c\in C\}$,
and therefore it's infinite. For $c\in C$, $(b_{c,j},b_{c,j+1},...,b_{c,m-1},b_{c,0},...,b_{c,j-1})$
and $w_j(g_0,...,g_{k_*-1})$ satisfy the same properties that $(b_{c,0},...,b_{c,m-1})$
and $w(g_0,...,g_{k_*-1})$ satisfy. Therefore, if $(l_1[c] : c\in C)$
is constantly $j>0$, then by conjugating and moving to $w_j(g_0,...,g_{k_*-1})$,
we may assume that $(l_1[c] : c\in C)$ is contantly $0$.

Let $l_2[c]<m$ be the maximal such that $\rho_{c,0}=...=\rho_{c,l_2[c]}$,
so wlog $l_2[c]=l_*$ for every $c\in C$. For $l<k_*$ let $\eta_{1,l}$
be $\bold{B}(f_l)$ if $l<k_{**}$ and $\nu_l$ if $l\in [k_{**},k_*)$
(we might also denote it by $\rho_l$ in this case). As $f_l \in H_3$
for $l<_{k_**}$, and $\rho_l \notin \{\bold{B}(f) : f\in H_3\}$
for $l\in [k_{**},k_*)$, it follows that $l_1<k_{**} \leq l_2<k_* \rightarrow \eta_{1,l_1} \neq \eta_{1,l_2}$.
Therefore, $(\eta_{1,l} : l<k_*)$ is without repetition. 

Now let $\eta_{2,l}$ be defined as follows:

1. If $l<k_{**}$, let $\eta_{2,l}$ be $\eta_2$ from definition
6(b)(2) for $f_l$ and $g_l$.

2. If $l\in [k_{**},k_*)$, let $\eta_{2,l}=\eta_{1,l}$.

Let $j(*)<\omega$ be such that:

a. $(\eta_{1,l} \restriction j(*) : l<k_*)$ is without repetition.

b. If $\eta_{1,l_1} \neq \eta_{2,l_2}$ then $\eta_{1,l_1} \restriction j(*) \neq \eta_{2,l_2} \restriction j(*)$
$(l_1,l_2<k_*)$.

c. If $\eta_{2,l_1} \neq \eta_{2,l_2}$ then $\eta_{2,l_1} \restriction j(*) \neq \eta_{2,l_2} \restriction j(*)$
$(l_1,l_2<k_*)$.

d. $j(*)>3m,k_*$.

e. $j(*)>n(l_1,l_2)$ for every $l_1<l_2<k_*$, where $n(l_1,l_2)$
is defined as follows:

1. If $k_{**} \leq l_1,l_2$, let $n(l_1,l_2)=0$.

2. If $l_1<k_{**}$ or $l_2<k_{**}$, let $(\nu_n^1 : n \in B_1)$
and $(\nu_n^2 : n\in B_2)$ be as in definition 6(b)(2) for $(f_{l_1},\eta_{1,l_1})$
and $(f_{l_2},\eta_{1,l_2})$, respectively. If there is no $\nu_n^1$
such that $\nu_n^1 \nleq \eta_{1,l_1}$ and no $\nu_n^2$ such that
$\nu_n^2 \nleq \eta_{1,l_2}$, let $n(l_1,l_2)=0$. Otherwise, there
is at most one $n\in B_1$ such that $\nu_n^1 \nleq \eta_{1,l_1}$
and $\nu_n^1 \leq \eta_{1,l_2}$ and there is at most one $m\in B_2$
such that $\nu_m^2 \nleq \eta_{1,l_2}$ and $\nu_m^2 \leq \eta_{1,l_1}$.
If there are $\nu_n^1$ and $\nu_m^2$ as above, let $n(l_1,l_2)=lg(\nu_n^1)+lg(\nu_m^2)+1$.
If there is $\nu_n^1$ as above but no $\nu_m^2$ as above, let $n(l_1,l_2)=lg(\nu_n^1)+1$,
and similarly for the dual case.

f. $j(*)>m(l_1,l_2)$ for every $l_1<l_2<k_{**}$ where $m(l_1,l_2)$
is defined as follows: Let $(\nu_n^1 : n\in B_1)$ and $(\nu_m^2 : m\in B_2)$
be as in definition 6(b)(2) for $(f_{l_1},\eta_{1,l_1})$ and $(f_{l_2},\eta_{1,l_2})$,
respectively. As $\eta_{1,l_1} \neq \eta_{1,l_2}$, $|\{\nu_n^1 : n\in B_1\} \cap \{\nu_m^2 : m\in B_2\}| <\aleph_0$,
let $s(l_1,l_2)$ be the supremum of the length of members in this
intersection and let $m(l_1,l_2):=s(l_1,l_2)+1$.

We may assume wlog that $lg(\rho_{c,l_1[c]})>j(*)$ for every $c\in C$.
We now consider two possible cases (wlog $TV((\rho_{c,l} : l<m)$
is constant$)$ is the same for all $c\in C$):

\textbf{Case I: For every $c\in C$, $(\rho_{c,l} : l<m)$ is not
constant.}

In this case, for each such $c\in C$, $l_2[c]<m-1$ and $b_{c,l_2[c]} \in u_{\rho_{c,0}}=...=u_{\rho_{c,l_2[c]}}$.
Now $(b_{c,l_2[c]},b_{c,l_2[c]+1}) \in g_{k(l_2[c])}^{i(l_2[c])}$,
and as $\rho_{c,l_2[c]} \neq \rho_{c,l_2[c]+1}$, necessarily $k(l_2[c])<k_{**}$.
By the definition of $l_1[c]$ and the fact that $\rho_{c,l_1[c]}=\rho_{c,l_2[c]}$,
necessarily $\rho_{c,l_2[c]}<_* \rho_{c,l_2[c]+1}$.

For each $l<m-1$, if $lg(\rho_{c,l})<lg(\rho_{c,l+1})$, then either
$g_{k(l)}$ or $g_{k(l)}^{-1}$ is as in definition 6(2)(b), so letting
$n=lg(\rho_{c,l})$, $(\rho_{c,l},\rho_{c,l+1})$ here correspond
to $(\eta_1 \restriction n,\nu_n)$ there, and there are $(a^l,b^l,c^l,d^l,e^l)=(a_c^l,b_c^l,c_c^l,d_c^l,e_c^l)$
in our case that correspond to $(a_n,b_n,c_n,d_n,e_n)$ in 6(2)(b).
In the rest of the proof we shall denote those sequences by $(a^l,b^l,c^l,d^l,e^l)$,
as the identity of the relevant $c\in C$ should be clear. In addition,
one of the following holds:

1. $i(l)=1$ and $g_{k(l)}^{i(l)}(a^l)=b^l$.

2. $i(l)=-1$ and $g_{k(l)}^{i(l)}(b^l)=a^l$.

Similarly, for $l<m-1$, if $lg(\rho_{c,l})>lg(\rho_{c,l+1})$ then
the above is true modulo the fact that now $(\rho_{c,l},\rho_{c,l+1})$
correspond to $(\nu_n,\eta_1 \restriction n)$ and one of the following
holds:

1. $i(l)=1$ and $g_{k(l)}^{i(l)}(b^l)=e^l$.

2. $i(l)=-1$ and $g_{k(l)}^{i(l)}(b^l)=a^l$.

Therefore, if $l=l_2[c]$ then $lg(\rho_{c,l})<lg(\rho_{c,l+1})$,
so the first option above holds, and therefore $\rho_{c,l}$ is an
initial segment of $\eta_{1,k(l)}$.

If $l=m-1$, then $lg(\rho_{c,m-1})>lg(\rho_{c,0})=lg(\rho_{c,m})$
and therefore $\rho_{c,0}$ is an initial segment of $\eta_{1,k(m-1)}$.
It follows that $\rho_{c,0}=\rho_{c,l_2[c]}$ is an initial segment
of $\eta_{1,k(l_2[c])} \cap \eta_{1,k(m-1)}$. Recalling that $lg(\rho_{c,0})=lg(\rho_{c,l_1[c]})>j(*)$
and that $(\eta_{1,l} \restriction j(*) : l<k_*)$ is without repetition,
it follows that $k(m-1)=k(l_2[c])$.

We shall now prove that if $l_2[c]<m-1$ then $l_2[c]=m-2$. Let $(a^{l_2[c]},b^{l_2[c]},...)$
be as above for $l=l_2[c]$, so $b_{c,l_2[c]+1}=b^{l_2[c]}$, and
as $k(l_2[c])=k(m-1)$, we get $b_{c,l_2[c]+1}=b^{l_2[c]}=b^{m-1}$.
In order to show that $l_2[c]=m-2$, it suffices to show that $b^{m-1}=b_{c,m-1}$
(as the sequence of the $b_{c,l}$s is without repetition), which
follows from the fact that $\rho_{c,0}<_*\rho_{c,m-1}$.

As we assume that the word $w$ is reduced, and as $k(m-2)=k(l_2[c])=k(m-1)$,
necessarily $i(m-2)=i(m-1)$. We may assume wlog that $i(m-2)=i(m-1)=1$
(the proof for $i(m-2)=i(m-1)=-1$ is similar, as we can replace $w$
by a conjugate of its inverse).

Let $w'=w'(g_0,...,g_{k_*-1}):=g_{k(m-3)}^{i(m-3)} \cdot \cdot \cdot g_{k(0)}^{i(0)}$,
by the above considerations and as $b^{m-1}=b_{c,m-1}$, it follows
that $e^{m-1}=g_{k(m-1)}(b^{m-1})=g_{k(m-1)}^{i(m-1)}(b^{m-1})=g_{k(m-1)}^{i(m-1)}(b_{c,m-1})=b_{c,0}$.
We also know that $g_{k(m-2)}^{i(m-2)}(w'(b_{c,0}))=g_{k(m-2)}(w'(b_{c,0}))$
is {}``higher'' than $w'(b_{c,0})$. Therefore, $g_{k(m-2)}^{i(m-2)}(w'(b_{c,0}))=g_{k(m-2)}(w'(b_{c,0}))=b^{m-2}=b^{m-1}$.
It also follows that $w'(b_{c,0})=a^{m-2}=a^{m-1}$. Therefore, $w'(e^{m-1})=a^{m-1}$.

We shall now prove that if $l<m-2$ then $b_{c,l+1}=(g_{\nu_{k(l)}}^*)^{i(l)}(b_{c,l})$.
Assume that for some $l<m-2$, $g_{k(l)}(b_{c,l}) \neq g_{\nu_{k(l)}}^*(b_{c,l})$
and we shall derive a contradiction. Let $(a^{m-2},b^{m-2},...)$
be as before for $g_{k(m-2)}$, so $(b_{c,m-2},b_{c,0})=(a^{m-2},e^{m-2})$.

Case I (a): $k(l)=k(m-2)=k(m-1)$. As the $b_{c,i}$s are without
repetition, if $0<l<m-2$, then $b_{c,l} \notin  \{b_{c,m-2},b_{c,0}\}=\{a^{m-2},e^{m-2}\}=\{a^l,e^l\}$,
and of course, $b_{c,l} \notin \{b^l,c^l,d^l\}$ (as it is a {}``lower''
element). Therefore, $g_{k(l)}(b_{c,l})=g_{\nu_{k(l)}}^*(b_{c,l})$,
a contradiction. If $l=0$, then $g_{k(m-1)}^{i(m-1)}(b_{c,m-1})=b_{c,0}$
and $(b_{c,m-1},b_{c,0})=(b^{m-1},e^{m-1})$. If $i(0)=-i(m-1)$,
then by conjugating $g_{k(0)}$, we get a shorter word with infinitely
many fixed points, contradicting our assumption on the minimality
of $m$.

If $i(0)=i(m-1)=1$, then we derive a contradiction as in the case
of $0<l$.

Case I (b): $k(l) \neq k(m-2)=k(m-1)$. In this case, we know that
$g_{k(l)}$ almost coincides with $g_{\nu_{k(l)}}^*$, with the exception
of at most $\{a^l,b^l,c^l,d^l,e^l\}$. Let $\rho_c:=\rho_{c,0}=...=\rho_{c,m-2}$,
then necessarily $\rho_c \leq \nu_{k(m-1)}=\eta_{1,k(m-1)}$ (as $g_{k(m-1)}$
moves $b_{c,m-1}$ to a lower $u_{\rho}$ (namely $u_{\rho_c}$),
$\rho_c$ plays the role of $\eta_1 \restriction n$ in Definition
6 for $g_{k(m-1)}$). By our assumption, $lg(\rho_c)>j(*)$ and $(\eta_{1,l} \restriction j(*) : l<k(*))$
is without repetition, therefore $\eta_{1,k(l)} \restriction lg(\rho_c) \neq \rho_c$,
so $\rho_c \nleq \eta_{1,k(l)}$. Therefore, when we consider $f_{k(l)}$
and $\eta_{1,k(l)}$ in definition $6(b)(2)$, then $\rho_c$ has
the form $\nu_n$ for some $n$. By the choice of $j(*)$, it's then
impossible to have $\rho_c \leq \eta_{1,k(m-1)}$, a contradiction.

Therefore, $a^{m-2}=w'(g_0,...,g_{k_*-1})(e^{m-2})=w'(g_0,...,g_{k_*-1})(b_{c,0})=w'(g_{\nu_0}^*,...,g_{\nu_{k_*-1}}^*)(b_{c,0})=w'(g_{\nu_0}^*,...,g_{\nu_{k_*-1}}^*)(e^{m-2})$.
In the notation of the claim and definition 6(b)(2), $F_1(f_{k(m-2)})(a^{m-2})=e^{m-2}$,
therefore, by composing with $w'$, we obtain a word composed of permutation
of $u_{\rho_{c,m-2}}$ (in the sense of claim 4(f)) that fixes $e^{m-2} \in u_{\rho_{c,m-2}}$,
therefore, $m-3=0$ (or else we get a contradiction by claim 4(f)).

It follows that $w(g_0,...,g_{k_*-1})=g_{k(m-2)}g_{k(m-2)}g_{k(0)}^{i(0)}$
and $g_{k(0)}^{i(0)}(e^{m-2})=a^{m-2}$. Now, obviously $\rho_{c,m-2} \leq \eta_{1,k(m-2)}$,
so $\rho_{c,m-2} \nleq \eta_{1,k(0)}=\nu_{k(0)}$. By the definition,
$g_{\nu_{k(0)}}^* \restriction u_{\rho_{c,m-2}}=f_{\rho_{c,m-2},\nu_{k(0)} \restriction lg(\rho_{c,m-2})} \neq f_{\rho_{c,m-2},\rho_{c,m-2}}$.
Also $F_1(f_{k(m-2)}) \restriction u_{\rho_{c,m-2}}=g^*_{\nu_{k(m-2)}} \restriction u_{\rho_{c,m-2}}=f_{\rho_{c,m-2},\rho_{c,m-2}}$.
Therefore, we get the following: $a^{m-2}=g_{k(0)}^{i(0)}(e^{m-2})=(g_{\nu_{k(0)}}^*)^{i(0)}(e^{m-2})=(f_{\rho_{c,m-2},\nu_{k(0)} \restriction lg(\rho_{c,m-2})})^{i(0)}(e^{m-2})$
and $e^{m-2}=F_1(f_{k(m-2)})(a^{m-2})=f_{\rho_{c,m-2},\rho_{c,m-2}}(a^{m-2})$.
In conclusion, we get a contrdiction to claim 4(f), as we have a short
non-trivial word that fixes $e^{m-2}$.

\textbf{Case II: $(\rho_{c,l} : l<m)$ is constant for every $c\in C$
(so $l_2[c]=m-1$). }Let $\rho_c:=\rho_{c,0}=...=\rho_{c,m-1}$. If
$g_{k(l)}^{i(l)}(b_{c,l})=(g_{\nu_{k(l)}}^*)^{i(l)}(b_{c,l})$ for
every $l<m$, then we get a contradiction to claim 4(f). Therefore,
for every $c\in C$, the set $v_c=\{l<m : g_{k(l)}^{i(l)}(b_{c,l}) \neq (g_{\nu_{k(l)}}^*)^{i(l)}(b_{c,l})\}$
is nonempty. Without loss of generality, $v_c$ doesn't depend on
$c$, and we shall denote it by $v$. For every $l\in v$, if $i(l)=1$
then $(b_{c,l},b_{c,l+1}) \in \{(a^l,b^l),(b^l,e^l),(c^l,d^l)\}$,
if $i(l)=-1$ then $(b_{c,l},b_{c,l+1}) \in \{(b^l,a^l),(e^l,b^l),(d^l,c^l)\}$.

We shall now prove that for some $k<k_{**}$, $k(l)=k$ for every
$l\in v$. Suppose not, then for some $l_1<l_2 \in v$, $k(l_1) \neq k(l_2)$.
By the choice of $j(*)$, each of the following options in impossible:
$\rho_{c,l_1} \leq \eta_{1,l_1} \wedge \rho_{c,l_2} \leq \eta_{1,l_2}$,
$\rho_{c,l_1} \leq \eta_{1,l_1} \wedge \rho_{c,l_2} \nleq \eta_{1,l_2}$,
$\rho_{c,l_1} \nleq \eta_{1,l_1} \wedge \rho_{c,l_2} \leq \eta_{1,l_2}$
or $\rho_{c,l_1} \nleq \eta_{1,l_1} \wedge \rho_{c,l_2} \nleq \eta_{1,l_2}$.
Therefore we get a contradiction. It follows that $\{k(l) : l\in v\}$
is singelton, and we shall denote its only member by $k<k_{**}$.

Note that if $l_1 \in v$, $l_2 \in v$ is the successor of $l_1$
in $v$, $l_1+1<l_2$ and $c\in C$ then $b_{c,l_1+1} \neq b_{c,l_2}$
(recall that $(b_{c,l} : l<m)$ is withut repetition). We shall now
arrive at a contradiction by examining the following three possible
cases (in the rest of the proof, we refer to $l(*)$ from Definition
5(A)(e) as {}``the distance between $a$ and $b$'', and similarly
for any pair of members from some $u_{\eta}$):

Case II (a): $g_k$ is as in definition 6(b)(2)(K)(a). In this case,
for every $l\in v$, the only possibilities for $(b_{c,l},b_{c,l+1})$
are either of the form $(c,d)$ or $(d,c)$ (and not both, as we don't
allow repetition). As the distance between $c$ and $d$ is at most
$2$, we get a word made of $f_{\rho,\nu}$s of length $\leq m+1$
that fixes $c$, contradicting claim 4(f).

Case II (b): $g_k$ is as in definition 6(b)(2)(K)(c). Pick $c\in C$
such that $lg(\rho_c)$ is alseo greater than $m+l_*$ where $l_*$
is as in definition 6(b)(2)(K)(c) for $g_k$. As the sequence $(b_{c,l} : l<m)$
is without repetition, necessarily $1\leq |v| \leq 3$.

If $|v|=3$, then necessarily the sequences $(a,b,e)$ or $(e,b,a)$
occur in $(b_{c,l} : l<m)$, as well as $(c,d)$ or $(d,c)$. As the
distance between $a$ and $e$ is $1$ and the distance between $c$
and $d$ is $\leq 2$, we get a contradiction as before.

Suppose that $|v|=2$. If the sequence $(a,b,e)$ appears in $(b_{c,l} : l<m)$,
we get a contradiction as above. If $(a,e)$ (or $(e,a)$) and $(c,d)$
(or $(d,c)$) appear, we also get a contradiction as above. If $(a,b)/(b,a)$
and $(c,d)/(d,c)$ appear, as the distance between $a$ and $b$ is
$l_*$, we get a word made of $f_{\rho,\nu}$s of length $\leq m+l(*)$
fixing $c$, a contradiction to claim 4(f). Finally, if $|v|=1$ we
get a contradiction similarly.

Case II (c): $g_k$ is as in definition 6(b)(2)(K)(b). As in the previous
case, where the only non-trivial difference is when either $|v| \in \{1,2\}$
and the sequence $(a,b)/(b,a)$ appears in $(b_{c,l} : l<m)$, but
not as a subsequence of $(a,b,e)/(e,b,a)$. If for some $c$ this
is not the case, then we finish as before, so suppose that it's the
case for every $c\in C$. As the distance between $c$ and $d$ is
$\leq 2$, suppose wlog that $|v|=1$, $k=k(m-1)$ (by conjugating)
and the sequence $(b_{c,l} : l<m)$ ends with $a$ and starts with
$b$ or vice versa. Therefore, every $c\in C$ is of the form $a^n$
or $b^n$ (where $n\in B$ and $B$ is as in definition 6(b)(2) for
$g_k$) and either $g_{k(0)}^{i(0)} \cdot \cdot \cdot g_{k_{m-2}}^{i(m-2)}(a^n)=b^n$
or $g_{k(0)}^{i(0)} \cdot \cdot \cdot g_{k_{m-2}}^{i(m-2)}(b^n)=a^n$,
so the distance between $a^n$ and $b^n$ is $\leq m-1$. As $C$
is infinite, the distance between $a^n$ and $b^n$ is $\leq m-1$
for infinitely many $n \in B$. This is a contradiction to the assumption
from definition 6(b)(2)(K)(b) that the distance between $a^n$ and
$b^n$ is increasing. 

This completes the proof of claim 8. $\square$

\textbf{Claim 9: }There exists a Borel function $\bold{B}_4: U^U \rightarrow U^U$
such that for every $f\in \bold{B}_1$, $\bold{B}_4(f) \in G_f$.

\textbf{Proof: }As in {[}HwSh:1089{]}, and we comment on the main
point in the proof of claim 7. $\square$

\textbf{Definition 10: }Let $G$ be the subgroup of $Sym(U)$ generated
by $\{\bold{B}_4(f) : f\in H_3\} \cup \{g_{\nu}^* : \nu \in 2^{\omega} \setminus \{\bold{B}(f) : f\in H_3\}\}$.

\textbf{Claim 11: }$G$ is a maximal cofinitary group.

\textbf{Proof: }$G$ is cofinitary by claim 8, so it's enough to prove
maximality. Assume towards contradiction that $H$ is a counterexample
and let $f_* \in H\setminus G$, so $\bold{B}_4(f_*) \in G$, and
we shall denote $f^*=\bold{B}_4(f_*)$.

\textbf{Case I: }$f_* \in H_3$. In this case, by Definition 6, $\{a_n : n\in B\} \subseteq eq(f_*,f^*):=X$
(see thee relevant notation in definition 6), hence it's infinite.
Therefore, $f_*^{-1}f^* \restriction X$ is the ientity, but $f_*^{-1}f^* \in H$
and $H$ is cofinitary, therefore $f_*^{-1}f^*=Id$ so $f_*=f^* \in G$,
a contradiction.

\textbf{Case II: }$f_* \notin H_3$. By the definition of $H_3$,
$\bold{B}_{2,6}'(f_*)=0$, so the sequences $\bold{B}_{2,1}(f_*)=(A_n=A_n(f_*) : n\in \bold{B}_{2,2}(f_*))$,
$\bar{\rho_*^n}=(\rho_{n,i} : i<l_*)=(\rho_{a,f_*(a),i} : i<l_*)$
and $\bar{i}=(i_{l} : l<l_*)=(i_{a,f(a),l} : l<l_*)$ $(n\in \bold{B}_{2,2}(f_*), a\in A_n)$
are well-defined, and for every $l<l_*$, $(\rho_{n,l} : n\in \bold{B}_{2,2}(f_*))$
is $\leq$-increasing, so $\nu_l:=\underset{n\in \bold{B}_{2,2}(f_*)}{\cup} \rho_{n,l} \in 2^{\omega}$
is well-defined. Let $g=(g_{\nu_0}^*)^{i_0} \cdot \cdot \cdot (g_{\nu_{l_*-1}}^*)^{i_{l_*-1}} \in G_1$(we
may assume that it's a reduced product). Let $w_1=\{l<l_* : (\exists f_l \in H_3)(\nu_l=\bold{B}(f_l))\}$
and $w_2=l_* \setminus w_1$. For $l<l_*$, define $g_l$ as follows:

1. If $l\in w_1$, let $g_l=\bold{B}_4(f_l)$.

2. If $l\in w_2$, let $g_l=g_{\nu_l}^*$.

Let $g'=g_0^{i(0)} \cdot \cdot \cdot g_{l_*-1}^{i(l_*-1)}$. By the
definition of $G$, $g_0,...,g_{l_*-1} \in G$, hence $g' \in G$.

Again by Definition 6, if $l\in w_1$ then $g_l=F_1(f_l)$ $mod$
$I_1$ and $g_l^{-1}=F_1(f_l)^{-1}$ $mod$ $I_1$. Now suppose that
$g(a) \neq g'(a)$, then there is a minimal $l<l_*$ such that $(g_{\nu_0}^*)^{i(0)} \cdot \cdot \cdot (g_{\nu_l}^*)^{i(l)}(a) \neq (g_0)^{i(0)} \cdot \cdot \cdot (g_l)^{i(l)}(a)$.
Let $v_l=dif(g_{\nu_l}^*,g_l)$, then $a\in (g_0^{i{0}} \cdot \cdot \cdot g_{l-1}^{i(l-1)})^{-1}(v_l)$.
In order to show that $(g_0^{i{0}} \cdot \cdot \cdot g_{l-1}^{i(l-1)})^{-1}(v_l) \in I_1$
it suffices to observe that for $i\in w_1$, functions of the form
$g_i,g_i^{-1}$ map elements of $I_1$ to elements of $I_1$, therefore
it follow that $g=g'$ $mod$ $I_1$. It suffices to show that $eq(f_*,g) \notin I_1$,
as it will then follow that $eq(f_*,g') \notin I_1$, so $f_*^{-1}g'=Id$
on an $I_1-$positive set, hence on an infinite set. As $f_*^{-1}g' \in H$
and $H$ is cofinitary, $f_*^{-1}g'=Id$, an therefore $f_*=g' \in G$,
a contradiction. 

So let $n\in \bold{B}_{2,2}(f_*)$ and $a\in A_n=A_n(f_*)$, and observe
that $f_*(a)=g(a)$. Indeed, by the definition of $\bold{B}_{2,1}(f_*)$,
for every such $a$, $f_*(a)=((f_{\eta_1 \restriction n \rho_{n,0}}^{i_0}) \cdot \cdot \cdot (f_{\eta_ \restriction n \rho_{n,l_*-1}}^{i_{l_*-1}}))(a)$
(where $\eta_1$ is as in the definition of $\bold{B}_{2,1}(f_*)$).
It's now easy to verify that the last expression equals $g(a)$. It's
also easy to verify that $\underset{ n\in \bold{B}_{2,2}(f_*)}{\cup}A_n \notin I_1$,
therefore we're done.

$\square$

\textbf{Claim 12: }$G$ is Borel.

\textbf{Proof: }It suffices to prove the following subclaim:

Subclaim: There exists a Borel function $\bold{B}_5$ with domain
$Sym(U)$ such that if $g\in G$ then $\bold{B}_5(g)=(g_0,g_1,...,g_m)$
such that $G\models "g=g_0^{i_0}g_1^{i_1} \cdot \cdot \cdot g_m^{i_m}"$
for some $(i_0,...,i_m) \in \{-1,1\}^{m+1}$. 

Proof: By the definition of $G$, if $g\in G$ then there are $m$,
$f_0,...,f_m \in \bold{A}_1$ (possibly with repetition) and $i_0,...,i_m \in \{-1,1\}$
such that $g=g_0^{i_0} \cdot \cdot \cdot g_m^{i_m}$ where each $g_i$
is either of the form $\bold{B}_4(f_i)$ for $f_i \in H_3$ (in this
case, let $\nu_i:=\bold{B}(f_i)$) or $g_{\nu_i}^*$ for $\nu_i \in 2^{\omega} \setminus \{\bold{B}(f) : f\in H_3\}$. 

Now if $n$ is greater than $m!$, then for some $u\subseteq 2^n$
such that $|u|\leq m!<\frac{2^n}{2}$, for every $\rho \in 2^n \setminus u$
we have:

a. For every $l\leq m$, $g_l \restriction u_{\rho}=f_{\rho,\nu_l \restriction lg(\rho)}$.

b. $g\restriction u_{\rho}$ can be represented as $f_{\rho,\nu_0 \restriction lg(\rho)}^{i_0} \cdot \cdot \cdot f_{\rho,\nu_m \restriction lg(\rho)}^{i_m} \in Sym(u_{\rho})$.

d. By claim 4(f), the above representation of $g\restriction u_{\rho}$
is unique.

Therefore, from $g$ we can Borel-compute $((\nu_i \restriction n : n<\omega) : i<m)$
hence $(\nu_i : i<m)$.

As $H_3$ is Borel and $\bold{B}$ is injective, the sets $\{\bold{B}(f) : f\in H_3\}$
and $2^{\omega} \setminus \{\bold{B}(f) : f\in H_3\}$ are Borel.
Now if $\nu_i \in 2^{\omega} \setminus \{\bold{B}(f) : f\in H_3\}$,
we can Borel compute $g_i=g_{\nu_i}^*$. If $\nu_i \in \{\bold{B}(f) : f\in H_3\}$,
then $\nu_i=\bold{B}(f_i)$ and we can Borel-compute $f_i$ (by applying
$\bold{B}_{-1}$ from definition 3(e)) hence $\bold{B}_4(f_i)$.

$\square$

\textbf{\large References}{\large \par}

{[}HwSh1089{]} Haim Horowitz and Saharon Shelah, A Borel maximal eventuallly
different family, arXiv:1605.07123.

{[}FFT{]} Vera Fischer, Sy David Friedman, Asger Toernquist, Definable
maximal cofinitary groups, arXiv:1603.02942.

{[}Ka{]} Bart Kastermans, The complexity of maximal cofinitary groups,
Proc. of the Amer. Math. Soc., Vol. 137(1), 2009, 307-316

{[}Ma{]} A. R. D. Mathias, Happy families, Ann. Math. Logic \textbf{12
}(1977), no. 1, 59-111. MR 0491197.

$\\$

(Haim Horowitz) Einstein Institute of Mathematics

Edmond J. Safra Campus,

The Hebrew University of Jerusalem

Givat Ram, Jerusalem, 91904, Israel.

E-mail address: haim.horowitz@mail.huji.ac.il

$\\$

(Saharon Shelah) Einstein Institute of Mathematics

Edmond J. Safra Campus,

The Hebrew University of Jerusalem

Givat Ram, Jerusalem, 91904, Israel.

Department of Mathematics

Hill Center - Busch Campus,

Rutgers, The State University of New Jersey.

110 Frelinghuysen road, Piscataway, NJ 08854-8019 USA

E-mail address: shelah@math.huji.ac.il
\end{document}